\magnification=1200
\input amssym.def
\input amssym.tex
 \font\newrm =cmr10 at 24pt
\def\bul{\raise .9pt\hbox{\newrm .\kern-.105em } }

 \def\fr{\frak}
 \font\sevenrm=cmr7 at 7pt
 
 \def\Z{{\Bbb Z}}

 \baselineskip 20pt
 
 \def\h{\hbox{ }}

 \def\u{{\fr u}}

 \def\m{{\fr m}}
 \def\n{{\fr n}}

 \def\ss{{\fr s}}
 
 \def\b{{\fr b}}
 
 \def\hh{{\fr h}}

 \def\g{{\fr g}}

 \def\q{{\fr q}}

 \def\<{\le}
 \def\>{\ge}

 \def\s{{\h\subset\h}}
 
 \def\vs{\vskip }

 \def\mapright#1
  {\smash{\mathop
  {\longrightarrow}
  \limits^{#1}}}

 \def\kk#1{{\kern .4 em} #1}
 \def\vs{\vskip 1pc}

\font\twelverm=cmr12 at 14pt
\hsize = 31pc
\vsize = 45pc
\overfullrule = 0pt
\font\sevenrm=cmr10 at 7pt
\font\smallfont=cmr10 at 9pt
\font\smallbold=cmb10 at 9pt

\font\authorfont=cmr10 at 9pt

\centerline{\twelverm Fomenko--Mischenko Theory, Hessenberg
Varieties,}
\centerline{\twelverm and Polarizations}
\vskip 1.5pc
\baselineskip=11pt
\vskip8pt
\centerline{\authorfont BERTRAM KOSTANT} 
\vskip 2pc
\baselineskip=14pt
\noindent{\smallbold Abstract.} {\smallfont The symmetric algebra $S(\g)$ over a
Lie algebra
$\g$ has the structure of a Poisson algebra. Assume $\g$ is 
 complex semisimple.  Then results of Fomenko--Mischenko (translation
of invariants) and A.~Tarasev construct a polynomial subalgebra
${\cal H} = \Bbb C[q_1,\ldots,q_b]$ of $S(\g)$ which is maximally
Poisson commutative. Here $b$ is the dimension of a Borel
subalgebra of $\g$. Let $G$ be the	adjoint group of $\g$ and let
$\ell = \hbox{\rm rank}\,\g$. Using the Killing form, identify $\g$ with its
dual so that any $G$-orbit $O$ in $\g$ has the structure (KKS) of
a symplectic manifold and $S(\g)$ can be identified with the affine
algebra of $\g$. 

An element $x\in \g$ will be called strongly regular if $\{(dq_i)_x\},\,i=1,\ldots,b$, are linearly
independent. Then the set $\g^{\hbox{{\sevenrm sreg}}}$ of all strongly regular elements is
Zariski open and dense in $\g$ and also $\g^{\hbox{{\sevenrm sreg}}}\s \g^{\hbox{{\sevenrm reg}}}$
where $\g^{\hbox{{\sevenrm reg}}}$ is the set of all regular elements in
$\g$. A Hessenberg variety is the $b$-dimensional
affine plane in $\g$, obtained by translating a Borel subalgebra by a
suitable principal nilpotent element. Such a variety was 
introduced in [K2]. Defining $\hbox{\rm Hess}$ to be a particular Hessenberg
variety, Tarasev has shown that $\hbox{\rm Hess}\s \g^{\hbox{{\sevenrm sreg}}}$. 

Let $R$ be the set of all regular $G$-orbits in $\g$. Thus if $O\in
R$, then $O$ is a symplectic manifold of dimension $2n$ where $n=
b-\ell$. For any $O\in R$ let $O^{\hbox{{\sevenrm sreg}}} = \g^{\hbox{{\sevenrm sreg}}}\cap O$. One
shows that $O^{\hbox{{\sevenrm sreg}}}$ is Zariski open and dense in $O$ so that
$O^{\hbox{{\sevenrm sreg}}}$ is again a symplectic manifold of dimension $2n$. For
any $O\in R$ let $\hbox{\rm Hess}(O) = \hbox{\rm Hess}\cap O$. One proves that 
 $\hbox{\rm Hess}(O)$
is a Lagrangian submanifold of $O^{\hbox{{\sevenrm sreg}}}$ and that $$\hbox{\rm Hess} =
\sqcup_{O\in R} \hbox{\rm Hess}(O).$$ The main result of this paper is to show
that there exists simultaneously over all $O\in R$,
an explicit polarization (i.e., a ``fibration" by Lagrangian
submanifolds) of $O^{\hbox{{\sevenrm sreg}}}$ which makes
$O^{\hbox{{\sevenrm sreg}}}$ simulate, in some sense, the cotangent bundle of
$\hbox{\rm Hess}(O)$. 
}

\vskip 6pt

\noindent  {\smallbold Keywords}: symplectic geometry, geometric quantization, Poisson manifolds, 
 symplectic manifolds, Lagrangian submanifolds, Poisson algebras, group actions, invariant theory, group actions
   on affine varieties, rings and algebras
\vskip 6pt
\noindent {\smallbold MSC classification codes}: 53D50, 53D05, 53D17, 53D12, 17B63, 17B66, 16Wxx

\vskip 1pc
\baselineskip 16pt
\centerline{\bf 0. Introduction}\vskip 1pc 
\rm
{\bf 0.1}. \rm  Let
$\g$ be a complex semisimple Lie algebra and put $\hbox{\rm rank}\,\g =
\ell$. Let
$G$ be the adjoint group of
$\g$. As one knows the symmetric algebra $S(\g)$ has the structure
of a Poisson algebra. Identify $\g$ with its dual using the Killing
form $(x,y)$ so that we can also regard $S(\g)$ as the algebra of
polynomial functions on $\g$. But then $\g$ inherits the structure
of a Poisson manifold. The corresponding symplectic leaves are the
adjoint orbits $O$ of $G$. For any $\varphi \in S(\g)$ let
$\xi_{\phi}$ be the``Hamiltonian" vector field on $\g$. If $x\in
\g$,  then one knows $$(\xi_{\varphi})\in T_x(O), $$ where $O$ is the
adjoint orbit of $x$. 

Let $\g = \n_- + \hh + \n$ be a standard triangular decomposition
of $\g$. Let $\b = \hh +\n$ (resp. $\b_- = \hh + \n_-$). Let $b=
\hbox{\rm dim}\,\b$ (resp. $n = \hbox{\rm dim}\,\n$) so that $\hbox{\rm dim}\,\g = b+n$. Let
$I_j\in S(\g)^G,\,j=1,\ldots,\ell, $ be homogeneous generators of
the algebra of $G$-invariants in $S(\g)$. Let $d_j$ be the degree
of $I_j$. One knows that  $$\sum_{j=1}^{\ell} d_j = b. \eqno (0.1)$$ For any
$u\in
\g$ let $\partial_y$ be the directional partial derivative on $\g$
defined by
$y$. Thus by (0.1) one obtains $b$ polynomials $\q_1,\ldots, q_b$
on $\g$ by considering all $(\partial_y)^k\,I_j$ where
$j=1,\ldots,\ell,$ and $k= 0,\ldots, d_j-1$. Let ${\cal H}_y$ be 
the subalgebra of $S(g)$ generated by the $q_i,i=1,\ldots,b$. Then

\vs {\bf Theorem 0.1}. (Fomenko--Mischenko) {\it ${\cal H}_y$ is 
Poisson commutative}.\vs
 
{\bf 0.2.} Let $\Pi =\{\alpha_1,\ldots,\alpha_{\ell}\}$ be the set
of simple positive (i.e., with respect to $\b$) roots and let
$\{e_{\alpha_i},\,i=1,\ldots,\ell,\}$ be corresponding roots
vectors.  Let $\{w,e,f\}$ be the S-triple	whose span $\u$ is the
principal TDS where $e = \sum_i^{\ell} e_{\alpha_i}$ and $w\in \hh$
is defined so that $\alpha_i(w) = 2$ for $i=1,\ldots,\ell$. 

Let $\g^{\hbox{{\sevenrm reg}}}$ be the dense Zariski open set of all regular elements
$x\in \g$ ($x$ is regular if $\hbox{\rm dim}\,\g^x =\ell$). Fix $y\in \hh$.
We introduce the following terminology. An element $z\in \g$ will be
said to be strongly regular if $\{(dq_i)_z,\,i=1,\ldots,b\}$ are
linearly independent. An old criterion of ours for regularity
implies $$\g^{\hbox{{\sevenrm sreg}}} \s \g^{\hbox{{\sevenrm reg}}}.\eqno (0.2)$$ It is immediately
obvious that if $\g^{\hbox{{\sevenrm sreg}}}$ is Zariski open in $\g$ and is Zariski
dense if $\g^{\hbox{{\sevenrm sreg}}}$ is not empty. However it is true but not
obvious that $\g^{\hbox{{\sevenrm sreg}}}$ is not empty. 

\vs {\bf Theorem 0.2.} (Fomenko--Mischenko)
{\it $\g^{\hbox{{\sevenrm sreg}}}$ is not empty. In fact if $e_1\in \Bbb C^{\times}$,
then $e_1\in \g^{\hbox{{\sevenrm sreg}}}$}. 

\vs The Fomenko--Mischenko proof of Theorem 0.2 involves a case-by-case argument over all complex simple Lie
algebras. In our paper here we give a case-independent proof (see
Theorem 2.11) using the representation theory of the principal TDS
$\u$.

One immediate consequence of Theorem 0.2 is that the polynomials
$q_i$ are algebraically independent so that $${\cal H}_y = \Bbb
C[q_1,\ldots,\q_b]$$ is a polynomial ring in $b$-variables. 

We will refer to a translate of a Borel subalgebra of the form $e_1
+ \b_-$ as a generalized Hessenberg variety. In the case at hand if
$N_-\s G$ is the subgroup corresponding to $\n_-$, then $e_1 + \b_-$
is  stable under the adjoint action of $N_-$. With positive and
negative roots reversed we studied this action in [K2]. One outcome
was the existence of a section of the action of $G$ on the set of
regular $G$-orbits. To fix matters here let $e_1\in \Bbb
C^{\times}\,e $ be normalized so that $(e_1,f) =1$ and put $$\hbox{\rm Hess}=
e_1 + \b_-\,.$$ \vskip .5pc In an all too brief note [T], A.A. Tarasev
proved that
${\cal H}_y$ was a maximal Poisson commutative subalgebra of
$S(\g)$. This solved the question of maximality raised by E.~Vinberg.
  In addition Tarasev generalized Theorem 0.2 by proving
that $$\hbox{\rm Hess}\s \g^{\hbox{{\sevenrm sreg}}}.\eqno (0.3)$$ 
Let $$\Phi:\g\to \Bbb C^b$$ be
the morphism defined by putting $\Phi(x) =
(q_1(x),\ldots,q_{b}(x))$. In [T] Tarasev implicitly proves 
that the restriction $$\Phi: \hbox{\rm Hess}\to \Bbb C^b\eqno (0.4)$$ is an
algebraic isomorphism. 

\vs {\bf 0.4.} In the present paper we 
apply the machinery above to simultaneously polarize a Zariski
dense open set in all maximal $G$-orbits $O$, and in doing so
simulate a cotangent bundle structure on these Zariski open sets.

Any $G$-adjoint orbit $O$ is a symplectic manifold, specifically
with respect to KKS structure in the complex category. One has
$\hbox{\rm dim}\,O \leq 2n$ and $O$ is an orbit of regular elements if and
only if $\hbox{\rm dim}\,O =2n$. Let $R$ be the set of all $G$-orbits of
regular elements. For any $O\in R$ let $$\hbox{\rm Hess}(O) = O\cap \hbox{\rm Hess}.$$ In
this paper we prove

  \vs {\bf Theorem 0.3.} {\it Let $O\in R$. Then
$\hbox{\rm Hess}(O)$ is a Lagrangian submanifold of $O$. In particular
$\hbox{\rm dim}\,\hbox{\rm Hess}(O) = n$. Furthermore one has the disjoint union
  $$\hbox{\rm Hess} =
\sqcup_{O\in R} \hbox{\rm Hess} (O).\eqno (0.5)$$ 
In addition $(0.5)$ is the
decomposition of $\hbox{\rm Hess}$ into $N_-$ orbits.}

\vs If $O$ is any
$G$-orbit,  let $O^{\hbox{{\sevenrm sreg}}} = O\cap \g^{\hbox{{\sevenrm sreg}}}$. It is immediate from
(0.2) that $O^{\hbox{{\sevenrm sreg}}}$ is empty if $O$ is not regular. On the
other  hand from Theorem 0.3 and (0.3) it follows that $O^{\hbox{{\sevenrm sreg}}}$
is not empty if $O\in R$. In such a case $O^{\hbox{{\sevenrm sreg}}}$ is necessarily
open and Zariski dense in $O$. In particular $O^{\hbox{{\sevenrm sreg}}}$ is then a
symplectic manifold of dimension $2n$ and $\hbox{\rm Hess}(O)$, by (0.3), is a
Lagrangian submanifold of $O^{\hbox{{\sevenrm sreg}}}$. 

Now for any $x\in \g^{\hbox{{\sevenrm sreg}}}$ let $Z_x\s T_x(\g)$ be the span
of
$(\xi_{q_i})_x,\,i=1,\ldots,b$. 

\vs {\bf Theorem 0.4.} {\it The
correspondence $x\mapsto \Z_x$ defines an $n$-dimensional involutive
distribution ${\cal Z}$ on $\g^{\hbox{{\sevenrm sreg}}}$, so that by Frobenius (in
the holomorphic category), one has a foliation $$\g^{\hbox{{\sevenrm sreg}}} =
\sqcup_{\lambda\in
\Lambda} L_{\lambda}\eqno (0.6)$$ for some parameter set
$\Lambda$, where the $n$-dimensional leaves 
$L_{\lambda}$ are maximal connected integral submanifolds of
${\cal Z}$. }
 On the other hand $\Phi$ defines a fibration
$\g^{\hbox{{\sevenrm sreg}}}$. Let $\Phi^{\hbox{{\sevenrm sreg}}} = \Phi|\g^{\hbox{{\sevenrm sreg}}}$ so that
$$\Phi^{\hbox{{\sevenrm sreg}}}:\g^{\hbox{{\sevenrm sreg}}}\to \Bbb C^b.\eqno (0.7)$$

\vskip .5pc For any
$x\in \hbox{\rm Hess}$ let $F_x$ be the fiber of $\Phi^{\hbox{{\sevenrm sreg}}}$ over $\Phi(x)$.
That is $$F_x = (\Phi^{\hbox{{\sevenrm sreg}}})^{-1}(\Phi(x))$$ so that $$\g^{\hbox{{\sevenrm sreg}}} =
\sqcup_{x\in \hbox{\rm Hess}} F_x. \eqno (0.8)$$ It is of course clear that
$F_x$ is a Zariski closed (not necessarily irreducible) subvariety
of
$\g^{\hbox{{\sevenrm sreg}}}$.  

For any $x\in \hbox{\rm Hess}$ let $$\Lambda_x = \{\lambda\in \Lambda\mid
L_{\lambda}\s F_x\}.$$ 
\vskip .5pc {\bf Theorem 0.5.} {\it
$L_{\lambda}$, for any $\lambda\in \Lambda$, is an irreducible,
Zariski closed, nonsingular, $n$-dimensional subvariety of
$\g^{\hbox{{\sevenrm sreg}}}$. Furthermore $$\Lambda = \sqcup_{x\in
\hbox{\rm Hess}}\Lambda_x.$$ In addition $\Lambda_x$, for any $x\in
\hbox{\rm Hess}$, is a finite set and one has  $$F_x = \sqcup_{\lambda\in
\Lambda_x} L_{\lambda}.\eqno (0.9)$$ Moreover $F_x$ is a nonsingular
$n$-dimensional Zariski closed subvariety of $\g^{\hbox{{\sevenrm sreg}}}$ and $(0.9)$
is both the decomposition of $F_x$ into the union of its irreducible
components and simultaneously the decomposition of $F_x$ into its
connected (with respect to both its Zariski and ordinary Hausdorff 
topology) components.}

\vskip 6pt Our final result is that the
maximal Poisson commutative subalgbra ${\cal H}_y$ of $S(\g)$ leads to
 a simultaneous polarization of $O^{\hbox{{\sevenrm sreg}}}$ for all regular
$G$-orbits $O$.

\vs {\bf Theorem 0.6.} {\it Let $O\in R$ and let
$x\in \hbox{\rm Hess}(O)$. Then $F_x\s O^{\hbox{{\sevenrm sreg}}}$. In fact $F_x$ is a
Lagrangian submanifold of $O^{\hbox{{\sevenrm sreg}}}$ and $$\g^{\hbox{{\sevenrm sreg}}} = \sqcup_{x\in
\hbox{\rm Hess}(O)} F_x,  \eqno (0.10)$$ thereby defining a
polarization of the $2n$-dimensional symplectic manifold $O^{\hbox{{\sevenrm sreg}}}$. 
Even more the Lagrangian submanifold $\hbox{\rm Hess}(O)$ of $O^{\hbox{{\sevenrm sreg}}}$ is
transversal to all the Lagrangian fibers $F_x,\,x\in \hbox{\rm Hess}(O)$, so
that $(0.10)$ simulates on $O^{\hbox{{\sevenrm sreg}}}$ the structure of the cotangent
bundle of
$\hbox{\rm Hess}(O)$.} 

\vskip 1.3pc \centerline{\bf Poisson structure and the
generalized Hessenberg variety}\vskip 8pt 
\centerline{\bf 1.
Poisson bracket on $\g$ and the principal TDS} \vskip .5pc {\bf
1.1.} Let
$\g$ be a complex semisimple Lie algebra and put $\hbox{\rm rank}\,\g =
\ell$. Let
$G$ be the adjoint group of
$\g$. Identify
$\g$ with its dual using the Killing form $(x,y)$. The symmetric algebra $S(\g)$ over $\g$ then
identifies with the algebra of polynomial functions on $\g$, where if $x,y\in \g$,  then $x(y) = (x,y)$.
If
$q$ is any holomorphic function on
$\g$ (e.g., elements of $S(\g)$) and $x\in \g$,  let $dq(x)\in \g$ be
defined so that for any $z\in \g$, $$(dq(x),z) = {d\over dt}(q(x+tz))_{t=0}.\eqno(1.1)$$

Let
$I_j,\,j=1,\ldots,\ell$, be homogeneous generators of $S(\g)^{G}$.
Let $d_j = m_j +1$ be the degree of $I_j$. The following old result of ours
is an immediate extension of Theorem 9, p.~382 in [K2].

\vskip .5pc {\bf Theorem 1.1.} {\it For any $x\in \g$ and $I\in S(\g)^G$ one
has $$dI(x)\in \hbox{\rm Cent}\,\g^x.\eqno (1.2)$$ Furthermore
$dI_j(x),\,j=1,\ldots,\ell$, is a basis of $\g^x$ if and only if
$x\in \g$ is regular. In particular if $x$ is regular semisimple, then $\g^x$ is
the unique Cartan subalgebra which contains $x$ so that $dI_j(x),\,j=1,\ldots,\ell$, is
a basis of the Cartan subalgebra $\g^x$. }

\vs {\bf Proof.} Let
$a\in G^x$ and
$z\in
\g$.  Then $$\eqalign{(a\cdot dI(x),z)&= (dI(x),a^{-1}\cdot z)\cr
&=  {d\over dt}(I(x+ta^{-1}\cdot z))_{t=0}\cr &= 
{d\over dt}(I(a^{-1}\cdot (x+t z))_{t=0}\cr &={d\over dt}
(I(x+t z))_{t=0}\,\,\hbox{since $I$ is  
$\hbox{\rm Ad}\,G$-invariant}\cr&=(dI(x),z)\cr}$$ so that $$dI(x)\in \g^{G^x}. \eqno
(1.3)$$ But $x\in \g^x$ so that $$dI(x)\in \g^x. \eqno (1.4)$$ But (1.3) and
(1.4) yield (1.2). 
      
But now by Theorem 9, p.~382 in [K2], one
has $dI_j(x),\,j=1,\ldots,\ell$, are linearly independent if and only if
$x$ is regular.  But, by definition, $x$ is regular if and only if $\hbox{\rm dim}\,\g^x = \ell$.
But, clearly, this proves the theorem. QED

\vs  {\bf 1.2.} Let $\hh$ be a Cartan subalgebra of
$\g$ and let $\Delta$ be the set of roots for $(\hh,\g)$. For each $\varphi\in \Delta$ let
$e_{\varphi}$ be a corresponding root vector. Let $\Delta_+\s \Delta $ be a choice of
positive roots and let $\Pi = \{\alpha_1,\ldots,\alpha_{\ell}\}$ be the
set of simple positive roots. Let $w\in \hh$ be the unique element such that
$\alpha(w) =2 $ for any $\alpha\in \Pi$. Let $f= \sum_{j=1}^{\ell}
e_{-\alpha_j}$ so that $$[w,f] = -2 f. \eqno (1.5)$$ One has that $w$ is regular
semisimple and $\g^w = \hh$. Let
$e\in
\sum_{j=1}^{\ell}\Bbb C\,e_{\alpha_j}$ be such that $\{w,e,f\}$ is a
principle $S$-triple spanning a principal TDS $\u$. Then by [K1] one has
a direct sum $$\g = \oplus_{j=1}^{\ell} \m_j, \eqno (1.6)$$ where
$$\hbox{\rm dim}\,\m_{j} = 2m_j +1\eqno (1.7)$$ is an $\hbox{\rm ad}\,\u$ irreducible module.

 Let
$\{z_j\},\,j=1,\ldots,\ell$, be a basis of $\hh$ such that $$\m_{j}\cap\hh
= \Bbb C z_j.\eqno (1.8)$$ Let $\b_-$ (resp.~$\b$) be the Borel subalgebra spanned by
$\hh$ and $\{e_{-\varphi} (\hbox{resp. $e_{\varphi}$}) \},\,\varphi\in \Delta_+$. Let $\n_- =
[\b_-,\b_-]$ (resp.~$\n = [\b,\b]$). Put $\m_{j,-} = \b_{-}\cap \m_j$. Let $$z_{j\,k} =
(\hbox{\rm ad}\,f/2)^k\,z_j,\,k=0,\ldots,m_j,\eqno (1.9)$$ so that (since $\n_-$ is the span of
 $ \hbox{\rm ad}\,w$ eigenvectors with negative eigenvalues) one has

\vs {\bf Proposition 1.2.} {\it The set $\{z_{j\,k}\},\,k=0,\ldots,m_j,$ is a basis of
$\m_{j,-}$ and $$\{z_{j\,k}\},\,\, j=1,\ldots,\ell,\,\,\,k= 0,\ldots,m_j \eqno
(1.10)$$ is a basis of $\b_-$.}

\vs For any $t\in \Bbb C$ let $u_t =
\hbox{\rm exp}\,t/2\,\,f$ so that $$u_t\cdot w = w + t\,f.\eqno (1.11)$$ Then one notes that,
for $j=1,\ldots,\ell,$
$$u_t\cdot z_j = z_{j\,0} + t\,z_{j\,1} +\cdots + t^{m_j}\,z_{j\,m_j}. \eqno (1.12)$$ But
then by the usual Vandermonde argument one has

 \vs {\bf Proposition 1.3.} {\it Let
$j=1,\ldots, \ell$. If
$c_i\in
\Bbb C, i=1,\ldots, d,$ are distinct numbers where $d\geq d_j$, then $$u_{c_i}\cdot
z_j,\,\,i=1,\ldots,d,\,\,\hbox{spans} \,\,\m_{j\,-}.\eqno (1.13)$$}

\vs Obviously $w + t\,f$
is regular semisimple, by (1.11), and the Cartan subalgebra $\g^{w+t\,f}$ equals 
$u_t\cdot \hh $ so that, by (1.12), $$\{z_{j\,0} + t\,z_{j\,1} +\cdots +
t^{m_j}\,z_{j\,m_j},\,\,j=1,\ldots,\ell\}\,\,\hbox{is a basis of the Cartan
subalgebra}\,\,\g^{w+t\,f}.\eqno (1.14)$$  
Another basis of $\g^{w+ t\,f}$ is given by
Theorem 1.1. That is, $$dI_j(w+tj),\,j=1,\ldots,\ell,\,\,\hbox{is also a basis of}\,\,
\g^{w+ t\,f}.\eqno (1.15)$$
 Now we recall that the Coxeter number $h$ of $\g$ is the maximal value of
$d_j,\,j=1,\ldots,\ell$. We have proved 

\vs {\bf Theorem 1.4.} {\it For any $t\in \Bbb C$
and $j=1,\ldots \ell$, one has $dI_j(w + tf)\in \b_-$. Furthermore $\b_-$ is spanned by
$dI_j(w + tf)$ for all
$j=1,\ldots,\ell$, and $t\in \Bbb C$. In fact it is already spanned by these elements
where $t$ is restricted to take on a finite set of values whose cardinality is greater
than or equal to the Coxeter number $h$.}

 \vs {\bf 1.3.} Let $\hbox{\rm Hol}(\g)$ be the algebra of
holomorphic functions on $\g$. Now if
$p,q\in \hbox{\rm Hol}(\g)$ one defines the Poisson bracket $[p,q]\in \hbox{\rm Hol}(\g)$ so that
for any
$x\in \g$,
$$[p,q](x) = (x,[dp(x),dq(x)]).\eqno (1.16)$$ One also defines the (holomorphic)
Hamiltonian vector field $\xi_p$ on $\g$ so that $$(\xi_p(q))(x) = [p,q](x).\eqno (1.17)$$
One notes that 
$$\eqalign{(\xi_p(q))(x)&= ([x,p(x)],dq(x))\cr &= {d\over
dt}(q(x+t[x,dp(x)])_{t=0}\cr}\eqno (1.18)$$ so that one has 

\vs {\bf Proposition 
1.5.} {\it Let $p\in \hbox{\rm Hol}(\g)$ and $x\in \g$. Then $(\xi_p)_x$ is tangent, at $x$, to
the adjoint orbit of $x$,  and in fact $(\xi_p)_x$ is the tangent vector, at $t=0$, to the
curve $(\hbox{\rm Ad}\,\hbox{\rm exp}\,-t\,dp(x))(x)$. That is, $$(\xi_p)_x = - [dp(x),x].\eqno (1.19)$$}
\indent {\bf
Remark 1.6.} Note that if
$p = z\in
\g$,  then for any $x\in \g$, $$dz(x) = z\eqno(1.20)$$ so that (1.18) becomes
$$(\xi_z(q))(x)=
 {d\over dt}(q(x+t[x,z]))_{t=0}. \eqno (1.21)$$ In particular $$(\xi_z)_x = -[z,x], \eqno
(1.22)$$ and hence of course $$T_x(O) = \{(\xi_z)_x\mid z\in \g\}.\eqno (1.23)$$ \vskip 1pc
\centerline{\bf 2. Coadjoint orbits and Fomenko--Mischenko Theory}\vskip 8pt {\bf 2.1.} The main
theorem (Theorem 2.11) of this section is due to Fomenko and Mischenko. Their proof is case-by-case
verification over all simple Lie algebras. Here we give a general proof using results in [K1] 
on the adjoint action of a principal TDS on $\g$ (see [K1]).

 Let
$x,z\in \g$. If the context leads to no confusion, we may identify $z$ with the tangent
vector
$(\partial_z)_x$ at $x$ where, for
$q	\in \hbox{\rm Hol}(\g)$ one has $$(\partial_z)_x q = {d\over dt }(q(x+ tz))_{t=0}.$$ 

\vskip .5pc Let $O$ be the adjoint (= coadjoint) orbit containing $x$. We recall that
$O$ has a symplectic structure (KKS) denoted by $(O,\omega_O)$, where if $\omega_x$ is the
value of $\omega_O$ at $x$, then for $y,z\in \g$, $$\omega_x(-[z,x],-[y,x]) =
(x,[y,z]).\eqno (2.1)$$ See \S 5.2, p. 180--183 in [K3] and
Theorem 5.3.1, p.~184 in [K3]. 

Let $\hbox{\rm Hol}(O)$ be the Poisson algebra of holomorphic functions on
$O$. If
$\varphi\in \hbox{\rm Hol}(O)$, the corresponding Hamiltonian vector field $\xi_{\varphi}$ on
$O$ is such that for any
$v\in T_x(O)$,
$$\omega_x((\xi_{\varphi })_x,v) = v\,\varphi.\eqno(2.2)$$ See (4.1.3), p.~166 in [K3]. 
Now by (1.22) we may choose $p\in \hbox{\rm Hol}(\g)$ so that $v = -[y,x]$ where $y = dp(x)$. Hence 
$$\omega_x((\xi_{\varphi})_x,-[y,x]) = (-[y,x]\varphi)(x).\eqno (2.3)$$ Now assume that
$\varphi = q|O$ where $q\in \hbox{\rm Hol}(\g)$. Then, by (1.18)
$$\eqalign{\omega_x((\xi_{\varphi})_x,-[y,x])&=  {d\over
dt}(q(x+t[x,dp(x)])_{t=0}\cr & = (\xi_p(q))(x)\cr &= [p,q](x)\cr &= (x,[dp(x),dq(x)])\cr 
&= \omega_x([-dq(x),x],v).\cr}\eqno(2.4)$$ Hence by the nonsingularity of $\omega_x$ and
(1.19) one has $$\eqalign{(\xi_{\varphi})_x &=[-dq(x),x]\cr &= (\xi_q)_x.\cr}\eqno (2.5)$$

\vskip .5pc As an immediate consequence of (2.5) one has

 \vs {\bf Proposition 2.1.} {\it
Let $p,q\in \hbox{\rm Hol}(\g)$ and let $O$ be an adjoint orbit. Then $$[p,q]\mid O =
[p\mid O,q\mid O].\eqno (2.6)$$ }\vskip .5pc Let $V\s S(\g)$ be a finite-dimensional space
of polynomial functions. For any $x\in \g$ let $$\g(V,x) = \{dp(x)\mid p\in V\}$$ so
that
$\g(V,x)$ is a subspace of $\g$. One notes that if $p_i,\,i=1,\ldots,\hbox{\rm dim}\,V$, is a basis $V$
and
$z_j,\,j=1,\ldots,\hbox{\rm dim}\,\g$, is a basis of $\g$, then $$\hbox{\rm dim}\,\g(V,x) =
  \hbox{\rm rank}\,M(V,x), \eqno
(2.7)$$ where $M(V,x)$ is the $\hbox{\rm dim}\,V\times \hbox{\rm dim}\,\g$ matrix $M_{i\,j}(V,x)$ given by 
$$M_{i\,j}(V,x) = (\partial_{z_j}p_i)(x).\eqno (2.8)$$ Now let $$m(V) = \max_{x\in\g}\,\,
\hbox{\rm dim}\,\g(V,x), \eqno (2.9)$$ and let $\g(V) = \{x\in \g\mid \hbox{\rm dim}\,\g(V,x) = m(V)\}$ so that
clearly 

\vs {\bf Proposition 2.2.} {\it Let $V\s S(\g)$ be any finite-dimensional space of
polynomial functions where $V\neq 0$. Then $\g(V)$ is a nonempty Zariski open subset of $\g$.}

 \vs
{\bf 2.2.} Let $n = \hbox{\rm card}\,\Delta_+$ and let $b = \ell + n$ so that 
 $$\eqalign{\hbox{\rm dim}\,\g &= \ell + 2\,n\cr
 \hbox{\rm dim}\,\b_- &=
b,\cr}\eqno (2.10)$$  and of course if $O$ is any adjoint orbit, then $$\hbox{\rm dim}\,O\leq
2\,n\eqno (2.11)$$ and $$\hbox{one has equality in (2.11) $\iff\,\, O$ is an orbit
of regular elements}.\eqno (2.12)$$ Let $\g^{\hbox{\sevenrm{reg}}}$ be the set of regular elements in $\g$
so that $\g^{\hbox{\sevenrm {reg}}}$ is a nonempty Zariski open set in $\g$.

 \vs 
{\bf
Proposition 2.3.} {\it Assume that
$V\s
S(\g)$ is a finite-dimensional space of Poisson commuting polynomial functions. Then
$$m(V)\leq b. \eqno (2.13)$$}
\indent  {\bf Proof.} Since the intersection of two
nonempty Zariski open sets in $\g$ is again a nonempty Zariski open set it suffices to
prove that $$b\geq \hbox{\rm dim}\,\g(V,x),\,\,\forall x\in \g^{\hbox{\sevenrm{reg}}}\eqno (2.14)$$
  Let $x\in \g^{\hbox{\sevenrm{reg}}}$ and
let $O$ be the (2n)-dimensional adjoint orbit containing $x$. Now considering tangent and
cotangent spaces for submanifolds. Let
$$\nu:T_x^*(\g)
\to T_x^*(O)\eqno (2.15)$$ be the surjection defined by the embedding of $O$ into $\g$. But
the kernel of $\nu$ is clearly $\ell$-dimensional. Hence the kernel of the restriction of
$\nu$ to $T_x^*[V] = \{dp_x\mid p\in V\} $ is at most $\ell$-dimensional. But, by (2.6),
$\nu(T_x^*[V])$ corresponds to an $\omega_x$-isotropic subspace of $T_x(O)$ under the
isomorphism $T^*_x(O)\to T_x(O)$ defined by $\omega_x$. Hence  $\hbox{\rm dim}\,\nu(T_x^*[V])\leq
r$. Thus $b \geq \hbox{\rm dim}\,T_x^*[V]$. But of course $\hbox{\rm dim}\,T_x^*[V] = \hbox{\rm dim}\,\g(V,x)$.
QED

\vs {\bf 2.3.} Let $j=1,\ldots,\ell, \,\, t\in \Bbb C$,  and $u\in \g$. Consider the
polynomial function on $\g$ whose value at $x\in \g$ is given by $I_j(t\,u + x)$. Note
that, over all $t\in \Bbb C$, one obtains a finite-dimensional subspace $V_{j,u}$ of
$S(\g)$. Indeed $V_{j,u}$ is spanned by the homogeneous polynomials $I_{j,u,k}(x),
k=0,\ldots, m_j,$ and constants, where we write $$I_j(t\,u + x) = I_j(t\,u) +
\sum_{k=0}^{m_j}\, I_{j,u,k}(x)\,t^k\eqno (2.16)$$ where $$\hbox{\rm deg}\,\,I_{j,u,k}(x) =
d_j-k.\eqno (2.17)$$ Now put $$V_u = \sum_{j=1}^{\ell}\,V_{j,\,u}.\eqno (2.18)$$

\vskip .5pc
{\bf Theorem 2.4} (Fomenko--Mischenko) {\it $V_u$ is Poisson commutative for any
$u\in \g$.} 

\vs {\bf Proof.} We must show that, for any $x\in \g$ and any $p,q\in V_u$,
$$ (x,[dp(x),dq(x)]) = 0.\eqno (2.19)$$ We first show that one has (2.19) if $p(x) =
I_j(t\,u +x)$ and $q = I_k(t\,u +x)$ where $j,k =1,\ldots,\ell$. Indeed, since exterior
differentiation commutes with translation one has  $ dI_k(t\,u +x)\in \hbox{\rm Cent}\,\g^{t\,u
+ x}$ for any $k =1,\ldots,\ell$, by Theorem 1.1. Thus $[dp(x),dq(x)] = 0$, establishing
(2.19) for this case. Now assume that $s,t\in \Bbb C$ are distinct. Let $p(x) = I_j(s\,u +
x)$. Then for any $z\in \g$ one has $$\eqalign{(s\,u + x, [dp(x),z]) &= 
([s\,u + x, dp(x)],z)\cr &= 0\cr}\eqno (2.20)$$ since $dp(x)\in \g^{s\,u +x}$ by
Theorem 1.1. But then using this argument twice one has that $[dp(x),dq(x)]$ is the Killing
form, orthogonal to both $s\,u + x$ and $t\,u +x$ if we put $q(x) = I_k(t\,u +x)$. But,
since $s\neq t$, $x$ is in the span of $s\,u + x$ and $t\,u +x$. This proves (2.19).
QED

\vs {\bf 2.4.} If $a\in G$ and $u\in \g$ it is clear that with respect to the adjoint
action of $a$ on $S(\g)$ one has $a\cdot V_u = V_{a\cdot u}$. It follows therefore that the
integer $m(V_u)$ depends only on the conjugacy class of $u$. We recall
$b = \ell + n$ so that $b$ is the dimension of a Borel subalgebra of $\g$. In
particular, recalling the notation in Theorem 1.4, one has $$b = \hbox{\rm dim}\,\b_- .\eqno (2.21)$$
By Proposition 2.3 one has $$m(V_u)\leq b\eqno (2.22)$$ for any $u\in \g$. Let $${\cal R}
= \{u\in \g\mid m(V_u)= b\}. \eqno (2.23)$$ 

\vskip .5pc {\bf Theorem 2.5.} {\it Any
principal nilpotent element of $\g$ lies in $\cal R$. }

\vs {\bf Proof.} By
conjugation it suffices to show that $f\in {\cal R}$ where $f$ is the principal nilpotent
element given in Theorem 1.4. But, by Theorem 1.4, one has $$\g(V_f,w) = \b_-,\eqno (2.24)$$
proving the theorem. QED

\vs {\bf Remark 2.6.} Assume that $u\in {\cal R}$. Then by
definition there exists $x\in \g$ and $q_{u,i}\in V_u,\,i=1,\ldots,b$, such that the
differentials $(dq_{u,i})_x$ are linearly independent. But this implies that the
polynomials
$q_{u,i}$ are algebraically independent. Thus if $A_u$ is the subalgebra of $S(\g)$
generated by the $q_{u,i}$,  it follows that $A_u$ is Poisson commutative and, as an
algebra, is given as the polynomial algebra $$A_u = \Bbb C[q_{u,1},\ldots,q_{u,b}]\eqno
(2.25)$$ in
$b$-variables.

\vs {\bf Theorem 2.7.} {\it ${\cal R}$ is a nonempty Zariski open subset of
$\g$. Furthermore
${\cal R}$ is closed under multiplication by $\Bbb C^{\times}$.} 

\vs {\bf Proof.} The
last statement is obvious from the definition of $V_u$. Now ${\cal R}$ is nonempty by 
Theorem 2.5. Let
$u\in {\cal R}$. Let $q_{u,i}$ and $x$ be as in Remark 2.6. Replacing $u$ by $v\in \g$ in
the definition of $q_{u,i}$ it is clear, from the matrix argument in (2.8), that the
set
$$\ss =
\{v\in
\g\mid (dq_{v,i})_x,\,i=1,\ldots,b,\,\,\hbox{are linearly independent}\}$$ is a Zariski
open neighborhood of $u$. But this proves the theorem. QED

\vs But now one has \vs {\bf
Theorem 2.8.} {\it ${\cal R}$ contains all regular semisimple elements.}

\vs {\bf Proof.}
Let $u$ be regular semisimple. By conjugacy we may assume that $u\in \hh$, using the
notation of \S 1.2. Recall (see \S 1.2) that $\n_- = [b_-,\b_-]$. Since ${\cal R}$ is
closed under scalar multiplication it suffices to show that $\lambda\,u\in {\cal R}$ for
some nonzero scalar $\lambda$. Let
$N_-\s G$ be the subgroup corresponding to $\n_-$. Consider the adjoint action of $N_-$ on
$\b_-$. Since 
$u$ centralizes no nonzero element in $\n_-$ one knows (e.g., by the Kostant--Rosenlicht
theorem, see e.g., bottom of p.~36 and 2.4.14 in [Sp]) that, for any $\lambda\in \Bbb C^{\times}$,
$$N_-\cdot
\lambda\,u =
\lambda\,u +
\n_-.\eqno (2.26)$$ In particular (using notation in \S 1.5) $\lambda\,u + f $ is conjugate
to
$\lambda\,u$. But any Zariski open neighborhood of $f$ contains $\lambda\,u + f$ for some
sufficiently ``small" $\lambda$. Hence $u\in {\cal R}$ by Theorem 2.7. QED

\vs {\bf 2.5.} 
Let $J_k\in S^k(\g)$. Let $x,y\in \g$. Then using the inner product on $S(\g)$ which
extends the Killing form one has $$\eqalign{J_k(x+ty) &= (J_k, (x+ty)^k/k!)\cr 
&=\sum_{j=0}^k\,(J_k, x^{k-j}/(k-j)!\,\,y^j/j!\,\,\,t^j/j!)\cr &=
\sum_{j=0}^k\,(1/j!)\,(\partial_y)^j\,J_k(x)\,t^j\cr &=\sum_{j=0}^k\,J_{k-j,y}(x)
t^j\cr}\eqno (2.27)$$ where
$$J_{k-j,y}= (1/j!)\,(\partial_y)^j\,J_k \in S^{k-j}(\g).\eqno (2.28)$$ Now as a function
of $x$ one has $dJ_k(x +ty)\in \g$ where for any $z\in \g $,
 $$\eqalign{(dJ_k(x+ty),z) &=
d/ds(J_k(x+ s\,z +t\,y))|_{s = 0}\cr&= \sum_{j=0}^k\,d/ds(J_{k-j,y}(x + s z))|_{s=0}\,t^j\cr
&= \sum_{j=0}^{k-1} (1/(k-j-1)!)((\partial_x)^{k-j-1}\,J_{k-j,y})(z)\cr}\eqno (2.29)$$ so
that $$\eqalign{dJ_k(x+ty) &= \sum_{j=0}^{k-1}
(1/(k-j-1)!)((\partial_x)^{k-j-1}\,J_{k-j,y})t^j\cr &=
\sum_{j=0}^{k-1}\,1/(j!\,(k-j-1)!)\,(\partial_y)^j(\partial_{x})^{k-j-1}J_k\,\,t^j .\cr}\eqno
(2.30)$$ 

Now for $j= 0,\ldots k-1$, write $v_{k-j-1} =
1/(j!\,(k-j-1)!)\,(\partial_y)^j(\partial_{x})^{k-j-1}J_n$ so that $$dJ_k(x+ty)=
\sum_{j=0}^{k-1} v_{k-j-1}t^j.\eqno (2.31)$$ Now assume that $J_k\in (S^k(\g))^G$ so
that
$[x+ty, dJ_k(x+ty)] = 0$ and hence equating coefficients of powers of $t$, one has

 \vs
{\bf Proposition 2.9.} (Fomenko--Mischenko) {\it $$\eqalign{[x, v_{k-1}] &= 0\cr
[v_0,y]&= 0,\cr}\eqno (2.32)$$ and for
$j=0,\ldots, k-2$, $$[x,v_{k-j-2}] = [v_{k-j-1},y].\eqno (2.33)$$}

\vskip 2pt Now fix $y$  to be a
regular element of $\hh$ and let $x= e$ recalling \S 1.2 so that $e$ is a principal
nilpotent element in the TDS $\u$. For $i =1,\ldots,h-1$, let $\b_j\s\b$ be the span of
all $e_{\varphi}$ where $(\varphi,w) = 2j$. Put $\b_0 = \hh$ so that $$\b
=\oplus_{j=0}^{h-1} \b_j.$$ One has, for $j=0,\ldots,h-2$,  $$\hbox{\rm ad}\,e:\b_{j}\to \b_{j+1}.
\eqno (2.34)$$ 
Now since $y \in \hh$ is regular it follows that $\b_j$ is stable under
$\hbox{\rm ad} \,y$ and $\hbox{\rm ad}\,y\mid \b_j$ is nonsingular for $j>0$. In particular 
  if $\hbox{\rm ad}_{\n}y=
\hbox{\rm ad}\,y\mid \n$,  then $\hbox{\rm ad}_{\n}y$ is invertible. 
Let $\zeta:\b\to\n $ be given by putting $\zeta=
-(\hbox{\rm ad}\,y)^{-1}\circ \hbox{\rm ad}\,e\mid\,\b$ so that for $i=0,\ldots,h-2$, 
$$\zeta:\b_i\to \b_{i+1},\eqno
(2.35)$$

 Now in the notation of (2.33) one notes that if $i=k-j-1,\,j=0,\ldots,k-1$, then
$$v_i\in
\m_i,\,i=0,\ldots,k-1,\eqno (2.36)$$ and (2.32) and (2.33) are the statements
$$\eqalign{[y,v_0]&=0\cr [e,v_{k-1}]&=0\cr \zeta(v_i) &= v_{i+1},\,i=0,\ldots,k-2.\cr}\eqno
(2.37)$$
\vskip .2pc {\bf Remark 2.10.} It is important to note that $\zeta$ is independent of
$k$ and $J_k$.

\vskip 6pt For $j=1,\ldots,\ell,$ and $i=0,\ldots,h-1, $ we define $v_i(I_j) = 0$ if
$i\geq d_j$ and $v_i(I_j) = v_i$ using the notation of (2.37) where $k=d_j$ and $J_k = I_j$.
One then has $$\{dI_j(e +t\,y)\mid t\in \Bbb C\} = \sum_{i=0}^{h-1} \Bbb C\,\, v_i(I_j)\eqno
(2.38)$$ and
$$\{dI_j(e +t\,y)\mid t\in \Bbb C\}\cap \m_i = \Bbb C v_i(I_j); \eqno (2.39)$$ and where
$v_h(I_j) = 0$, $$\zeta(v_i(I_j)) = v_{i+1}(I_j).\eqno (2.40)$$ 

 Recalling \S 2.1 one has
$$\g(V_y,e) = \oplus_{i=0}^{h-1} (\g(V_y,e))_i\eqno (2.41) $$ where $$(\g(V_y,e))_i=
\g(V_y,e)\cap \m_i\eqno (2.42)$$ and $(\g(V_y,e))_i$ is given by $$(\g(V_y,e))_i =
\sum_{i=1}^{\ell}\Bbb C v_i(I_j).\eqno (2.43)$$ 

 We can now prove the following result of
Fomenko--Mischenko (see, Lemma 4.3, p.~383 and Lemma 44,
p.384 in [F-M]). The proof of this result in [F-M] depends on the fact that (2.34)  (and hence
(2.35)) is surjective. The authors assert that this can be  proved by considering the question 
  case-by-case. We will give a general proof using the representation theory of the TDS $\u$.
Namely, one has that (2.34) is surjective since the spectrum of $\hbox{\rm ad}\,w$ on $\n$ is strictly
positive. 

\vs {\bf Theorem 2.11} (Fomenko--Mischenko).  {\it Let $y\in \hh$ be regular. Then
$$\g(V_y,e) = \b.\eqno (2.44)$$}
\indent  {\bf Proof.} The proof will be by induction on $i$
using (2.43). One has $(\g(V_y,e))_0 = \g^y$ by Theorem 1.1. But $\g^y = \hh$ and $\hh =
\b_0$. Assume inductively that $\sum_{m=0}^j\,\b_m\s \g(V_y,e)$ for $j\leq n-2$. Let
$v_{j+1}\in \b_{j+1}$ be arbitrary. By the surjectivity of (2.35) there exists $v_j\in
\b_{j}$ such that $\zeta(v_j) = v_{j+1}$. But by induction $v_j\in \g(V_y,e)_j$. Thus there
exists constants $c_k\in \Bbb C,\,k=1,\ldots,\ell$, such that $$v_j =\sum_{k=1}^{\ell}
c_k\,v_j(I_k),$$ But then $$v_{j+1} = \sum_{k=1}^{\ell}\,c_k\, v_{j+1}(I_k)$$ by (2.40).
Thus $v_{j+1}\in \g(V_y,e)$. QED 

\vs {\bf Remark 2.12. } Note that upon conjugating (2.44) by
an element in $\hbox{\rm exp}\,\hh$ we may replace 
$e$ in (2.44) by any element $e_1$ of the form $$e_1 = \sum_{i=1}^{\ell} b_i\,e_{\alpha_i}, $$
where all $b_i$ are in $\Bbb C^{\times}$. 

\vskip 1.5pc
\centerline{\bf 3. The generalized Hessenberg variety}\vskip 8pt
{\bf 3.1}  Let
$\q =
\Bbb C f +
\b$. We refer to an affine plane of the form 
$f + \b$ as a generalized Hessenberg variety and $\q$ as its linearization. We will also consider the
opposed linearized Hessenberg variety 
$\q_- =
\Bbb C e +
\b_-$. In essence the results in \S 3 are due to A.A.Tarasev. They are either implicit or
explicit in the very brief note [T].  For what we believe is greater clarity we will reestablish
Tarasev's results and place them in a context which will lead to the results of \S 4. The proof here is
along the lines leading to our result in [K2] that $f + \g^e$ is a section of the adjoint action of $G$
on
$\g^{reg}$. The generalized Hessenberg variety was introduced in [K2]. See \S 4 in [K2].

Clearly the
$b +1$-dimensional subspaces
$\q$ and
$\q_-$ are nonsingularly paired by the Killing form. Let $\q_-^{\perp}$ be the Killing form
orthocomplement of $\q_-$ in $\g$ so that $\q_-^{\perp}\s \n$ and $$\g = \q \oplus
\q_-^{\perp}. \eqno (3.1) $$ If $X$ is an affine variety,  then $A(X)$ will denote the affine algebra
of $X$. Consider $A(\q_-)$. By restricting the polynomial functions on $\g$ (namely $S(\g)$) to
$\q_-$ one has an exact sequence
$$0\longrightarrow(\q_-^{\perp})\longrightarrow S(\g)\longrightarrow A(\q_-) \longrightarrow
0 \eqno (3.2)$$ where $(\q_-^{\perp})$ is the ideal in $S(\g)$ defined by $\q_-^{\perp}$.
On the other hand one has the direct sum $$S(\g) = (\q_-^{\perp}) \oplus S(\q) \eqno
(3.3)$$ so that the restriction of the third map in (3.2) to
$S(\q)$ defines an algebra isomorphism
$$S(\q)
\to A(\q_-).\eqno (3.4)$$ Let $$Q:S(\g) \to S(\q)\eqno (3.5)$$ be the projection
defined by (3.3) so that for any $p\in S(\g)$ the
$$\hbox{image of both $p$ and $Q(p)$ in $A(\q_-)$ are the same.}\eqno (3.6)$$ One notes then
that
$$\eqalign{Q(e_{-\alpha_i}) &= c_i\,f,\,i=1,\ldots,\ell,\,\,\hbox{for some $c_i\in \Bbb
C^{\times}$}\cr Q([\n_-,\n_-])&= 0\cr Q& = \hbox{Id on $S(\b)$}.\cr}$$ The
decomposition (3.3) is clearly stable under $ad\,w$ so that $Q$ commutes with $ad\,w$. In
particular $Q$ maps $S(\g)^{w}$ into $S(\q)^{w}$. But of course $Q$ is the identity on
$S(\q)^{w}$ so that $$Q(S(\g)^{w}) = S(\q)^w.\eqno (3.7)$$

\vskip .5pc {\bf Remark 3.1.}
Since (3.3) is clearly a decomposition of graded vector spaces note that (3.7) on homogeneous
components may be written $$Q(S^m(\g)^{w}) = (S^m(\q)^w).$$

\vskip 6pt  Let
$e_1\s
\Bbb C e
$ be the normalization so that $$(e_1,f) = 1, \eqno (3.8)$$ and let $\hbox{\rm Hess}\s \q_-$ be the  fixed affine
variety defined by putting $$\hbox{\rm Hess} =
e_1 +\b_-.\eqno (3.9)$$ In particular $\hbox{\rm Hess}$ is a $b$-dimensional affine plane and $A(\hbox{\rm Hess})$
  is the 
affine ring of
$\hbox{\rm Hess}$. Again  restriction of functions defines a surjection $$\sigma_{\hbox{\sevenrm Hess}}:S(\g)\to
A(\hbox{\rm Hess}).\eqno (3.10)$$ If
$\sigma:S(\g) \to A(\q_-)$ is defined by restriction of functions (so that $\sigma$ is the
composite of $Q$ and the isomorphism (3.4)), then of course $$\sigma_{\hbox{\sevenrm Hess}} = 
\tau_{\hbox{\sevenrm Hess}}\circ \sigma,
$$ where $$\tau_{\hbox{\sevenrm Hess}}:A(\q_-) \to A(\hbox{\rm Hess})\eqno (3.11)$$ is defined by restriction of
functions.  We note that $$\sigma_{\hbox{\sevenrm Hess}}:S(\b) \to A(\hbox{\rm Hess})\eqno (3.12)$$ is an algebra
isomorphism of algebras $$\sigma_{\hbox{\sevenrm Hess}} (f) = 1.\eqno (3.13)$$

Now for any $k\in \Bbb Z_+$ let $S(\g)_{[k]}$ be the graded subspace of $S(\g)$ (with
homogeneous components  $S^m(\g)_{[k]}$) defined by putting $$S(\g)_{[k]} =\{p\in S(\g)\mid
\hbox{\rm ad}\,w (p) = k\,p\}.$$ 
But now if $q\in S(\q)^w$, clearly we may uniquely write, as a finite
sum,
$$q =
\sum_{k=0}\,f^k\,q_{[k]},\eqno (3.14)$$ where $$q_{[k]} \in S(\b)_{[k]}.\eqno (3.15)$$

 One
notes that if $q\in S^m(\q)^w$, then the sum in (3.14) can be taken for $k\leq m$, and one
has 
$$q_{[k]}\in S^{m-k}(\b)_{[k]}.\eqno (3.16)$$

 But now if $m>0$, one has $q_{[m]} = 0$ since 
of course $S^0(\b)_{[m]} = 0$ so that for $m>0$ and $q\in S^m(\q)^w$, one has $$q =
\sum_{k=0}^{m-1} f^k\,q_{[k]}\,\,\,\hbox{with}\,\, q_{[k]}\in  S^{m-k}(\b)_{[k]}.\eqno (3.17)$$

 Since the affine space $\hbox{\rm Hess}$ is a translation of $\b_-$ the tangent space to 
$\hbox{\rm Hess}$ at $e_1$ identifies with $\b_-$. Consequently, using the Killing form
nonsingular pairing of $\b$ and $\b_-$, one has an identification $$T^*(\hbox{\rm Hess}) = \b.\eqno
(3.18)$$ 

\vskip 6pt {\bf Proposition 3.2.} {\it Assume
$m>0$. Let
$p\in S^m(\g)^w$ and let
$q= Q(p)$ so that
$q\in S^m(\q)^w$. See Remark 3.1. 
Then 
$$\eqalign{\sigma_{\hbox{\sevenrm Hess}}(p) &=\sigma_{\hbox{\sevenrm Hess}}(q)\cr &=
\sigma_{\hbox{\sevenrm Hess}}\bigl(\sum_{k=0}^{m-1}q_{[k]}\bigr),\cr}\eqno (3.19)$$ where $q_{[k]}\in  S^{m-k}(\b)_{[k]}$.
In particular
$$q_{[m-1]}\in \b_{m-1} \eqno (3.20)$$ and in fact, using (3.18),  $$q_{[m-1]} =
d(\sigma_{\hbox{\sevenrm Hess}}q)_{e_1}.
 \eqno (3.21)$$ 
In addition if
$m>1$, then one has
$$\sum_{k=0}^{m-2}q_{[k]}\in S\bigl(\sum_{k=0}^{m-2}\b_k\bigr).\eqno (3.22)$$}
\indent {\bf Proof.} The
equality (3.19) follows from (3.6) and (3.13). Of course (3.20) is given by (3.17). If $z\in
\b_-$ then, since $\sigma_{\hbox{\sevenrm Hess}}(S^k(\b))$ vanishes at $e_1$ for all positive integers $k$,
clearly (3.18) implies
$$d/dt(q(e_1 + t\,z))_{t=0} = (q_{[m-1]},z).\eqno (3.23)$$ But this establishes (3.21). Finally write
$q_{[k]},\,k\leq m-2$, as a linear combination of the basis of
$S(\b)$ formed by all product monomials of the root vectors $e_{\varphi},\,\varphi\in \Delta_+$,
and a basis of  
$\hh$. It follows from the condition that $q_{[k]}\in S(\b)_{[k]}$ that the coefficient of every
monomial containing any $e_{\varphi}\in \sum_{j \geq m-1}\,\b_{j}$ is zero. But this implies
(3.22). 

\hfill QED

\vs {\bf 3.2.} Let $y\in \hh$ be regular. By Theorem 2.11 and Remark 2.12 one has
$$\g(V_y,e_1) = \b.\eqno (3.24)$$ In particular $$\hbox{\rm dim}\,V_y \geq b. \eqno (3.25)$$ However in the
notation of \S 1.1, Proposition 1.2, one has $$b = \sum_{j=1}^{\ell} d_j.\eqno (3.26)$$ But by
(2.18) one has $$V_y =\sum_{j=1}^{\ell} V_{j,y},\eqno (3.27)$$ where $$\hbox{\rm dim}\,V_{j,y}\leq d_j\eqno
(3.28)$$ by definition of $V_{j,\,y}$ in \S 2.3,  recalling that $V_{j,y}$ is spanned by the
homogeneous polynomials $I_{j,y,k},\,k=0,\ldots,m_j$, and $$I_{j,y,k}\in S^{d_j-k}.\eqno (3.29)$$ 
Thus one has

\vs {\bf Theorem 3.3.} {\it $$\eqalign{\hbox{\rm dim}\,V_y &=b\cr \hbox{\rm dim}\,V_{j,y} &= d_j. \cr}\eqno (3.30)$$
  Also
(3.27) is a direct sum and the homogeneous polynomials 
$I_{j,y,k},\,j=1,\ldots,\ell,\,\,\,k=0,\ldots,m_j,$ are a basis of $V_y$.}\vs

 We assume the
$d_j$ are nondecreasing in $j$ so that $d_{\ell}= h$ where we recall $h$ is the Coxeter number. Let the  
partition of $b$, dual to (3.26), be given as $$b = \sum_{m=1}^h r_m, \eqno (3.31)$$ where the $r_m$
are nonincreasing. It follows easily from (3.29) that $$r_m = \hbox{\rm dim}\,V_y^m, \eqno (3.32)$$ where 
$$V_y^m = V_y\cap S^m(\g). \eqno (3.33)$$ On the other hand from the representation theory of the
TDS $\u$ (yielding the surjectivity of (2.34)),  one readily has that $$r_m = \hbox{\rm dim}\,\b_{m-1}, \eqno
(3.34)$$  and hence proving

\vs {\bf Theorem 3.4.} {\it One has $$\hbox{\rm dim}\,V_y^m = \hbox{\rm dim}\,\b_{m-1},\,\,m = 1,\ldots,h.\eqno
(3.35)$$}

{\bf Remark 3.5.} For $m=1,\ldots,h$, consider the set of all pairs
$\{d_j,k\}$ in (3.29) such that $d_j-k = m$. For any such pair let $i = \ell + 1 -j$ and put
$$ J_{m,y,i} = I_{j,y,k}.\eqno (3.36)$$ Then one notes that
$$\{J_{m,y,i}\},\,i=1,\ldots,r_m,\,\,\,\hbox{is a basis of $V_y^m$}.\eqno (3.37)$$

 \vs {\bf 3.3.} As
above let
$y\in
\hh$ be regular. Let
${\cal H}_y\s S(\g)$ be the (Poisson commutative--see Theorem 2.4) subalgebra generated by the  
$b$-dimensional subspace $V_y$. In dealing with (3.36) it is convenient to simply order the pairs
$(m,i)$. Let ${\cal B} =\{1,\ldots,b\}$ and let ${\cal P} = \{(m,i)\mid m=1,\ldots,h-1,\,\,i =
1,\ldots ,r_m\}$. Recalling (3.31) let $${\cal B}\to {\cal P},\,\,\,\beta\mapsto
(m(\beta),i(\beta))\eqno (3.38)$$ be a bijection where if $\beta <\beta'$, then $ m(\beta)\leq
m(\beta') $ so that if $J_{y;\beta} = J_{m(\beta),y,i(\beta)}$, then (see Theorem 3.3)
$$\{J_{y;\beta}\mid \beta\in {\cal B}\}\,\,{\hbox {\rm is a basis of}}\,\, V_y \,\,{\hbox{\rm and also generates}
\,\, {\cal
H}_y}.\eqno (3.39)$$

\vskip 4pt {\bf Theorem 3.6.} {\it The function restriction map (see (3.9)) $${\cal
H}_y\to A(\hbox{\rm Hess})\eqno (3.40)$$ is an algebra isomorphism. Furthermore ${\cal H}_y$ is a
polynomial algebra. In fact $${\cal H}_y = \Bbb C[J_{y;1},\ldots,J_{y;b}].\eqno (3.41)$$
 Moreover
not only are the $J_{y;\beta}$ algebraically independent but in fact the differentials 
$$\{(dJ_{y;\beta})_{v}\},\,\,\beta\in {\cal B}\}\,\,\hbox{is a basis of}\,\, 
T^*_v({\hbox{\sevenrm Hess}})$$
$${\hbox{at any point}\,\, v\,\, \hbox{of the Hessenberg Hess--(not just at}\,\,
e_1).} \eqno(3.42)$$
 In fact,  even stronger,
(3.42) remains true if the $J_{y;\beta}$ are replaced by the restrictions $J_{y;\beta}\mid
\hbox{\rm Hess}$. Indeed the restrictions $J_{y;\beta}\mid
\hbox{\rm Hess},\,\beta\in {\cal B}$ define a ``coordinate system" on $\hbox{\rm Hess}$. In fact the map
$$\hbox{\rm Hess}\to \Bbb C^{b},\,\,\,v  \mapsto (J_{y;1}(v),\ldots,J_{y;b}(v))\eqno (3.43)$$ is an
algebraic isomorphism.}

\vs {\bf Proof.} Let $\beta\in {\cal B}$. For notational convenience put $p=
J_{y;\beta}$ so that $p = J_{m,y,i}$ where $m = m(\beta)$ and $i= i(\beta)$. Now recall the
notation of Proposition 3.2. Of course
$p\in S^m(\g)^w$.  As in Proposition 3.2 let $q = Q(p)$ so that $q\in S^m(\q)^w$.
 Write $q = q_{\beta}$ and $q_{[k]} = q_{\beta,[k]}$ so that by (3.18),$$\sigma_{\hbox{\sevenrm Hess}}(J_{y;\beta})
= \sigma_{\hbox{\sevenrm Hess}}(\sum_{k=0}^{m(\beta)-1}\,q_{\beta,[k]}).\eqno (3.44)$$ Also, by Proposition 3.2. one
has  $$q_{\beta, [m(\beta)-1]}\in \b_{m(\beta)-1}, \eqno (3.45)$$ and if
$m(\beta) > 1$, $$\sum_{k=0}^{m(\beta)-2}\,q_{\beta,[k]})\in
S\bigl(\sum_{k=0}^{m(\beta)-2)}\,\b_{k}\bigr).\eqno (3.46)$$ In adddition, by (3.18) and (3.21), $$q_{\beta,
[m(\beta)-1]} = d(\sigma_{\hbox{\sevenrm Hess}}(J_{y;\beta}))_{e_1}.\eqno (3.47)$$

\vskip .5pc For notational
convenience put $z_{\beta} = q_{\beta, [m(\beta)-1]}$. We now assert that $$\{z_{\beta}\mid
\beta\in {\cal B}\}\,\,\hbox{is a basis of $\b$}\eqno (3.48)$$ and hence in particular, by (3.45),
for $m=1\ldots,h$, $$\{z_{\beta}\mid m(\beta) = m\}\,\,\,\hbox{ is a basis of $\b_{m-1}$}.\eqno
(3.49)$$ Indeed, by dimension, Theorem 2.11 and (3.39), $$\{dJ_{y;\beta}(e_1)\mid \beta\in {\cal
B}\}\,\,\hbox{is a basis of $\b$}.\eqno (3.50)$$ However since $dJ_{y;\beta}(e_1)$ is an element of
$\b$ it is immediate from definitions that $$dJ_{y;\beta}(e_1) =
(d\sigma_{\hbox{\sevenrm Hess}}(J_{y;\beta}))_{e_1}.\eqno (3.51)$$ But this proves (3.48) (and (3.49)). Note also that
(3.50) and (3.51) establish that $\{\sigma_{\hbox{\sevenrm Hess}}(J_{y;\beta})\mid \beta\in {\cal B}\}$ (and a
fortiori $\{J_{y;\beta})\mid \beta\in {\cal B}\}$) are algebraically independent. 
 But now (3.48) implies
that
$S(\b)$ as a polynomial algebra can be given as
$$S(\b) =
\Bbb C[z_{1},\ldots,z_{b}].\eqno (3.52)$$
 On the other hand since $\hbox{\rm Hess}$ is the $e_1$
translate of $\b_-$ one has an algebra isomorphism $$\sigma_{\hbox{\sevenrm Hess}}:S(\b)\to A(\hbox{\rm Hess})\eqno
(3.53)$$ and hence the map $$\hbox{\rm Hess} \to \Bbb C^b,\,\,\,v\mapsto
(z_{1}(v),\ldots,z_{b}(v))\eqno (3.54)$$  is an algebraic isomorphism. But we wish to
show that the $z_{\beta}$ in (3.54) can be replaced by the $J_{y;\beta}$. To do this we first
show, inductively, that for all $\beta\in {\cal B}$ one has $$\sigma_{\hbox{\sevenrm Hess}}(\Bbb
C[z_1,\ldots,z_{\beta}]) \s \sigma_{\hbox{\sevenrm Hess}}(\Bbb C[J_{y;1},\ldots,J_{y;\beta}]).\eqno (3.55)$$  But
now if $m(\beta) = 1$, then (3.55) is true since if $i\leq \beta$, then $m(i) = 1$ and hence by (3.44) 
one has
$$\sigma_{\hbox{\sevenrm Hess}} (J_{y;i}) = \sigma_{\hbox{\sevenrm Hess}}(z_{i}).\eqno (3.56)$$

Now assume $m(\beta)>1$ (which
implies that $\beta >1$), and assume that (3.55) is true for $\beta-1$. But then by the induction
assumption and (3.46), one has 
  $$\sigma_{\hbox{\sevenrm Hess}}\bigl(\sum_{n=0}^{m(\beta)-2}\,q_{\beta,[n]}\bigr)\in
\sigma_{\hbox{\sevenrm Hess}}(\Bbb C[J_{y;1},\ldots,J_{y;m(\beta)-1}).\eqno (3.57)$$
  Here we are implicitly using
the obvious fact that if $k\leq m(\beta)-1$ and $m(\beta') = k$, then $\beta'\leq \beta-1$. But now, by
(3.44), (3.46) and (3.57) we may write $$\sigma_{\hbox{\sevenrm Hess}}(J_{y;\beta})= \sigma_{\hbox{\sevenrm Hess}}
 (z_{\beta}) + u, \eqno
(3.58)$$ where $u\in \sigma_{\hbox{\sevenrm Hess}}(\Bbb C[J_{y;1},\ldots,J_{y;m(\beta)-1}])$. Thus we may solve for
$\sigma_{\hbox{\sevenrm Hess}}(z_{\beta})$, establishing the inductive step (3.55). But then the remaining statements of
Theorem 3.6 are immediate consequences of (3.52), (3.53) and (3.54). QED

\vs By (3.29), (3.36), (3.39),
(3.41) and also (3.31), (3.32)  we can write down the Poincar\'e series $P_y(t)$ of ${\cal H}_y$. Namely
one has (clearly independent of the regular element $y\in\hh$)
$$\eqalign{P_y(t) &= 
\prod_{j=1}^{d_j} {1\over 1-t^1}\cdots {1\over 1-t^{d_j}}\cr &=\prod_{m=1}^h\,\,{1\over
(1-t^m)^{r_m}}.
\cr}\eqno (3.59)$$

\vskip 1pc \centerline{\bf 4. Strong
regularity and Zariski dense cotangent structure on regular orbits}\vskip 6pt
 {\bf 4.1.}
 Recall $n = b-\ell$ so that
$$b+n = \hbox{\rm dim}\,\g.\eqno (4.1)$$ Let ${\cal N} = \{1,\ldots,n\}$. Fix
   a regular semisimple element
$y\in\hh$. By (2.16) one has $$I_{j,y,0} =
I_j,\,\,j=1,\ldots,\ell, \eqno (4.2)$$ and hence by (3.36)
$$J_{d_j,y,\ell +1 -j} = I_j.\eqno (4.3)$$ Thus, recalling the
notation of Theorem 3.6, let ${\cal I}$ be the subset of
cardinality $\ell$ given by $${\cal I} = \{\beta \in {\cal
B}\mid m(\beta) = d_j,\,i(\beta) = \ell
+i-j,\,\,\,j=1,\ldots,\ell\}\eqno (4.4)$$ so that $${\cal B} = {\cal N}\sqcup {\cal I}$$ and 
$$\{J_{y;\beta},\,\,\beta\in {\cal I}\} =
\{I_j,\,j=1,\ldots,\ell \}.\eqno (4.5)$$ For notational
convenience let  $\{q_i,\, i\in {\cal B}\}$ be a reordering of
the $J_{y,\beta}$ so that one retains $$\{q_{\beta},\, \beta \in {\cal I}\} = \{I_j,\,\,j=1,\ldots,
\ell\}.\eqno (4.6)$$ Then, first of all let
$\Phi$ be the morphism of $\g$ to $\Bbb C^{b}$ given by $$\Phi(z) =
(q_1(z),\ldots q_{b}(z)).\eqno (4.7) $$ Then Theorem 3.6 asserts \vs
{\bf Theorem 4.1.} {\it The morphism $\Phi$ is surjective and in
fact the restriction $$\Phi: \hbox{\rm Hess}\to \Bbb C^b\eqno
(4.8)$$ is an algebraic isomorphism. In particular
the isomorphism $$(\Phi\mid \hbox{\rm Hess})^{-1}: \Bbb C^b\to
\hbox{\rm Hess}\eqno (4.9)$$ is a cross-section of $\Phi$}.

\vs Let
$$\g^{\hbox{\sevenrm sreg}} = \{ z\in \g\mid (dq_j)_z,\,j\in {\cal B}\,\,\hbox{be 
linearly independent}\}.\eqno (4.10)$$ Then $\g^{\hbox{\sevenrm sreg}}$ is not empty by
  Theorem 3.6, and in fact Theorem 3.6 asserts that  $$\hbox{\rm Hess}\s
\g^{\hbox{\sevenrm sreg}}.\eqno (4.11)$$ The elements in $\g^{\hbox{\sevenrm sreg}}$ are regular
by (4.6) and our criterion for regularity (see reference at
 the end of \S 1.1). Thus $$\g^{\hbox{\sevenrm sreg}}\s \g^{\hbox{\sevenrm sreg}}.\eqno (4.12)$$  We will refer
to the elements in $\g^{\hbox{\sevenrm sreg}}$ as strongly regular. 

Now since any $I\in S(\g)^G$  Poisson commutes with any $p\in S(\g)$
one has $\xi_{q_i} = 0$ for $i\in {\cal I}$. For $i \in {\cal N}$ let $\xi_i =
\xi_{\q_i}$. By Poisson commutativity $$[\xi_i,\xi_j] = 0,\,\,\,\forall i,j\in {\cal N}.\eqno (4.13)$$ On
the other hand if
$x\in
\g$ and
$O$ is the
$G$-adjoint orbit containing $x$ then, for $i\in {\cal N}$, 
$$(\xi_i)_x \in T_x(O)\eqno (4.14)$$ by (2.5). If $x\in \g^{\hbox{\sevenrm sreg}}$,  then
as noted in (2.12),  one has $$\hbox{\rm dim}\, O = 2n.\eqno (4.15)$$

 \vskip .3pc Let $R$ be the set of all 
 $G$-orbits in $\g^{\hbox{\sevenrm sreg}}$. If $x$ is strongly regular,  let $Z_x$ be the span
of $(\xi_i)_x,\,\,i\in {\cal N}$ so that $Z_x\s T_x(O)$. 

\vs {\bf Theorem 4.2.} {\it Let $x\in
\g^{\hbox{\sevenrm sreg}}$ and let $O\in R$ be the regular orbit containing $x$. Then $Z_x$ is a Lagrangian subspace
of $T_x(O)$ and $$\{(\xi_i)_x\},\,i\in {\cal N},\,\hbox{is a basis of $Z_x$.}\eqno (4.16)$$}
{\bf Proof.} Since 
$\{(dq_j)_x\},\,j\in {\cal B}$ are linearly independent and since $\{q_k\},\,k\in {\cal I}$,
are constant on $O$ it is immediate that $\{(dq_k)_x\},\,k\in {\cal I},$ are a basis of the
orthocomplement of $T_x(O)$ in $T^*_x(\g)$. But then necessarily the differentials
$\{d(q_i|O)_x\},\,i\in {\cal N},$ are a basis of $T^*_x(O)$. But this implies (4.16). But then $Z_x$
is Lagrangian by (4.13). QED

\vs {\bf 4.2.} For any $O\in R$ and $j\in \{1,\ldots,\ell\}$ let $I_j(O)$
be the constant value that $I_j$ takes on $O$. We recall (see Theorem 2, p. 360 in
[K2]) that if 
$\eta:R\to\Bbb C^{\ell}$ is the map given by putting $\eta(O) = (I_1(O),\ldots,I_{\ell}(O))$, then
$$\eta:R\to \Bbb C^{\ell}\,\,\,\,\hbox{is a bijection}.\eqno (4.17)$$ Now for any $O\in R$ let 
$$\hbox{\rm Hess}(O) = O\cap \hbox{\rm Hess}.$$ Elsewhere we have proved (reversing the roles of $\n$ and $\n_-$)

 \vskip 6pt
{\bf Theorem 4.3.} {\it Let $O\in R$. Then 

\indent\indent (a)  $\hbox{\rm Hess}(O)$ is the subvariety of $\hbox{\rm Hess}$ given by

                     \centerline{$\hbox{\rm Hess}(O) = \{ v\in \hbox{\rm Hess}\mid I_j(v) = 
     I_j(O),\,\,\,j=1,\ldots,\ell\}$;}
\indent\indent (b)  $\hbox{\rm Hess}(O)$ is a principal (i.e., with trivial isotropy subgroup) $N_-$ orbit in 
 \hbox{\rm Hess}
so that in particular $\hbox{\rm Hess}(O)$ is a nonsingular $n$-dimensional subvariety of $\hbox{\rm Hess}$.

\indent\indent (c) One has $$\hbox{\rm Hess} = \sqcup_{O\in R} \hbox{\rm Hess}(O)\eqno (4.18)$$ so that (4.18) is the
$N_-$-orbit decomposition of $\hbox{\rm Hess}$.}

\vs {\bf Proof.} Theorem 4.3 is an immediate consequence of
Theorem 7, p.~381, and Theorem 8, p.~382 in [K2] and Theorem 1.2. p.~109 in [K4]. QED

\vs 
Let $O\in R$. As an adjoint orbit we recall that $O$ is a $2n$-dimensional symplectic manifold. On the
other hand $\hbox{\rm Hess}(O)$ is an $n$-dimnesional submanifold of $O$ by (b) in Theorem 4.3. 
In fact \vs 

{\bf Theorem. 4.4.} {\it Let $O\in R$. Then $\hbox{\rm Hess}(O)$ is a Lagrangian submanifold of $O$.}

\vs {\bf Proof}. Let $v\in \hbox{\rm Hess}(O)$. Let $\tau_1, \tau_2\in T_v(\hbox{\rm Hess}(O)$. Then if, as in (2.2),
$\omega_v$ is the symplectic form $\omega_O$ at $v$, we are to show that $\omega_v(\tau_1,\tau_2)= 0$.
But by Theorem 4.3 there exists $z_i\in \n_-$, for $i=1,2$, such that $-[z_i,v] = \tau_i$. But then,
by (2.1), $$\omega_v(\tau_1,\tau_2) = (v, [z_2,z_1]).\eqno (4.19)$$ But $[z_2,z_1]\in [\n_-,\n_-]$ and
clearly $$ [\n_-,\n_-]\s \sum_{k=2}^{h-1}\g_{[-k]},\eqno (4.20)$$ where $\g_{[k]} = \{x\in\g\mid [w,x]
= 2k\,x\}$. On the other hand $v\in \g_{[1]} + \hh + \sum_{k=1}^{h-1}\g_{[-k]}$. Hence 
$(v,[z_2,z_1])=0$. Thus $\hbox{\rm Hess}(O)$ is Lagrangian. QED

\vs {\bf 4.2.} Now clearly $\g^{\hbox{\sevenrm{sreg}}}$ is a
nonempty (by e.g., (4.11)) Zariski open subset of $\g$. In particular $\g^{\hbox{\sevenrm{sreg}}}$ is a quasi-affine
nonsingular irreducible algebraic variety. Let
$\Phi^{\hbox{\sevenrm{sreg}}}= \Phi|\g^{\hbox{\sevenrm{sreg}}}$ (see (4.7)), so that by (4.8) one has the surjective
morphism $$\Phi^{\hbox{\sevenrm{sreg}}}:\g^{\hbox{\sevenrm{sreg}}} \to \Bbb C^b$$  and  
 $$\Phi^{\hbox{\sevenrm{sreg}}}:\hbox{\rm Hess}\to \Bbb C^b$$ is an
algebraic isomorphism. For any $c\in \Bbb C^{b}$ let $$F_c = (\Phi^{\hbox{\sevenrm{sreg}}})^{-1}(c)$$ so that $F_c$
is a closed subvariety of $\g^{\hbox{\sevenrm{sreg}}},$ (noting that variety in our notation here does not require
irreducibility) and $$\g^{\hbox{\sevenrm{sreg}}} = \sqcup_{c\in \Bbb C^b}\,F_c.\eqno (4.21)$$ 
Of course if
$c=(c_1,\ldots,c_b)\in \Bbb C^{b}$, then $$F_c = \{z\in \g^{sreg}\mid q_i(z) =
c_i,\,\,i=1,\ldots,b\}.\eqno (4.22)$$ On the other hand one knows that for any $z\in \g^{\hbox{\sevenrm{sreg}}}$ the
differentials $(dq_i)_z,\, i=1,\ldots,b,$ are linearly independent. 

\vs {\bf Theorem 4.5.} {\it $F_c$
is a nonsingular variety of dimension $n$ for any $c\in  \Bbb C^{b}$. }

\vs {\bf Proof.} Let $z\in
F_c$. Then Theorem 4.5 is an immediate consequence of Theorem 4,
\S 4 in Chapter III, p. 172 in [M] where $U\s \g^{\hbox{\sevenrm{sreg}}}$ is an affine neighborhood of $z$ and
$f_1,\ldots,f_{b}$ are the images of $q_1,\ldots,q_b$ in $A(U)$. QED

\vs  {\bf Theorem 4.6.} {\it Let
$c\in \Bbb C^b$. Then the analytic space $F_c$ is a nonsingular analytic manifold of dimension n.}

\vs
{\bf Proof.} This is immediate from Corollary 2 , \S 4 in Chapter III, p. 168 in [M]. 

\hfill QED 

 \vskip 3pt For
any $G$-orbit $O$ in $\g$ let $$O^{\hbox{\sevenrm{sreg}}} = \g^{\hbox{\sevenrm{sreg}}}\cap O.$$

\vskip .5pc {\bf Proposition 4.7.} {\it
Let $O$ be a $G$-orbit in $\g$. Then $O^{\hbox{\sevenrm{sreg}}}$ is nonempty if and only if $O\in R$. In fact if $O\in
R$, then $O^{\hbox{\sevenrm{sreg}}}$ is an open Zariski dense subvariety of $O$. In particular $O^{\hbox{\sevenrm{sreg}}}$
  is a $2n$-dimensional symplectic submanifold of $O$, and one has $$\g^{\hbox{\sevenrm{sreg}}} =
  \sqcup_{O\in R}\,O^{\hbox{\sevenrm{sreg}}}.\eqno
(4.23). $$}
 {\bf Proof.} The proposition is immediate from (4.11), (4.12), Theorem 4.3 and of
course from the fact that $\g^{\hbox{\sevenrm{sreg}}}$ is Zariski open in $\g$. QED

\vskip 3pt Recalling (4.6), let $j(\beta)
\in \{1,\ldots,\ell\}$ be defined for $\beta\in {\cal I}$ so that $$q_{\beta}= I_{j(\beta)}.\eqno
 (4.24)$$ Then for $O\in R$, let $$\Bbb C^b(O)= \{c=(c_1,\ldots,c_b)\in \Bbb C^b\mid 
c_{\beta} = I_{j(\beta)}(O)\,\,\,\forall \beta\in {\cal I}\}.$$
 Of course $$\Bbb C^b = \sqcup_{O\in
R} \Bbb C^b(O). \eqno (4.25)$$ Then one has the following fibration (with $n$-dimensional fibers)
of $O^{\hbox{\sevenrm{sreg}}}$ for any $O\in R$. 

\vs {\bf Theorem 4.8.} {\it Let $O\in R$. Then $$ O^{\hbox{\sevenrm{sreg}}} =
\sqcup_{c\in
\Bbb C^b(O)}\,F_c.\eqno (4.26)$$}
\indent {\bf Proof.} This is immediate from Theorem 2, p.~360 in [K2].
This result asserts that any element $x\in \g^{\hbox{\sevenrm{sreg}}}$ is uniquely determined, up to $G$-conjugacy, by
the vector $(I_1(x),\ldots,I_{\ell}(x))\in \Bbb C^{\ell}$, and any such vector can be achieved by some
$x\in \g^{\hbox{\sevenrm{reg}}}$. QED

 \vs Theorem 4.2 asserts that $x\mapsto Z_x$ for $x\in \g^{sreg}$ is an
$n$-dimensional distribution (in the sense of differential geometry) ${\cal Z}$ on the
analytic manifold $\g^{\hbox{\sevenrm{sreg}}}$. But then (4.13) asserts that ${\cal Z}$ is involutory. Thus by the
Frobenius theorem, in the complex analytic category, one has a foliation of $\g^{\hbox{\sevenrm{sreg}}}$ by a family 
${\cal L}$ of  maximal integral connected ($n$-dimensional) manifolds of ${\cal Z}$. For the validity of
the use of the Frobenius theorem in the complex analytic category, see Theorem 1.3.6, p.~30 in [V] and
the comment at the end of \S 1.3 in Chapter 1, p. 31, in [V]. We refer to the elements
$L$ of 
${\cal L}$ as leaves of ${\cal Z}$. Let $\Lambda$ be an index set for ${\cal L}$ so that ${\cal L} =
\{L_{\lambda}\mid \lambda\in \Lambda\}$, and one has $$\g^{\hbox{\sevenrm{sreg}}} =\sqcup_{\lambda\in \Lambda}
L_{\lambda}. \eqno (4.27)$$ Recalling the notation of Theorem 4.2 note that, by definition of integral
manifold, for any
$\lambda\in
\Lambda$ and $x\in L_{\lambda}$,  one has $$T_x(L_{\lambda}) = Z_x.\eqno (4.28)$$
 A complex algebraic variety has, besides the
Zariski topology, the ordinary Hausdorff topology, which following Chapter 7 in [S], we refer to as the
complex topology.

\vs {\bf Theorem 4.9.} {\it Let $\lambda\in \Lambda$. Then there exists a unique $c\in \Bbb
C^{b}$ such that $L_{\lambda}$ is an open, in the complex topology, submanifold of the fiber
$F_c$. In particular if
$x\in L_{\lambda}$, then $$T_x(F_c) = Z_x.\eqno (4.29)$$}
\indent {\bf Proof.} If $i,j = 1,\ldots,b$, then
by Theorem 2.4 one has $[q_i,q_j] = 0$. Thus $\xi_i\,q_j=0$. But then $q_j$ is constant on
$L_{\lambda}$ for $j=1,\ldots,b$. But this implies that there exists $c\in \Bbb C^b$ such that
$L_{\lambda} \s F_c$. Since both $L_{\lambda}$ and $F_c$ are (nonsingular) analytic manifolds of
dimension $n$, this implies that $L_{\lambda}$ is open in $F_c$ in the complex topology and also that
(4.28) implies (4.29). QED 

\vs Now for any $c\in \Bbb C^b$ let $$\Lambda_c = \{\lambda\in \Lambda\mid
\L_{\lambda}\s F_c\}.$$ Recalling (4.27) the following statement is an immediate corollary of 
Theorem 4.9. 

\vs {\bf Theorem 4.10.} {\it Let $c\in \Bbb C^{b}$. Then $$ F_c = \sqcup_{\lambda\in
\Lambda}L_{\lambda}.\eqno (4.30)$$ Moreover (4.30) is the decomposition of the fiber $F_c$ into its
connected  components with respect to its complex topology.} 

\vskip 4pt But now recall (see Theorem
4.5) that the Fiber $F_c$ is a nonsingular algebraic variety. Using the Zariski toplogy this leads to a
much more interesting statement than Theorem 4.10.

 \vs {\bf Theorem 4.11.} {\it
The leaf $L_{\lambda}$ is a (Zariski) closed nonsingular algebraic irreducible subvariety of
$\g^{\hbox{\sevenrm{sreg}}}$ for any $\lambda\in \Lambda$. Furthermore if $c\in \Bbb C^{b}$, then $\Lambda_c$ is finite
and the  decomposition (4.30) is both the decomposition of the fiber $F_c$ into the union of its
(algebraic) Zariski irreducible components and also the decomposition of $F_c$ into the union of
its Zariski connected components.}

\vskip 4pt {\bf Proof.} Let $c\in \Bbb C^b$ and for some index set
$\Gamma$ let
$\{F_c^{\gamma},\,\gamma\in \Gamma\}$ be the set of all Zariski connected components of $F_c$.  Thus
$$F_c=\sqcup_{\gamma\in \Gamma}F_c^{\gamma}.\eqno (4.31)$$ But since $F_c$ is nonsingular the set of
Zariski connected components of $F_c$ is the same as the set of (Zariski) irreducible components of
$F_c$. See Corollary 17.2, p. 74 in [B]. Thus $\Gamma$ is finite and hence all the $F_c^{\gamma}$ are
Zariski open and closed in $F_c$. But Zariski open implies complex open. Thus if
$\lambda\in
\Lambda_c$ and $\gamma\in \Gamma$ one has, since $L_{\lambda}$ is complex connected, either
$L_{\lambda}\s F_c^{\gamma}$ or
$L_{\lambda}\cap F_c^{\gamma}$ is empty, Thus there exists a subset $\Lambda_c^{\gamma}\s \Lambda_c$
such that
$$F_c^{\gamma} = \sqcup_{\lambda\in \Lambda_{c}^{\gamma}} L_{\lambda}.\eqno (4.32)$$ But Theorem 4.30
implies that (4.32) is the decomposition of $F_c^{\lambda}$ into its complex connected components. But
$F_c^{\gamma}$ is complex connected by Theorem, \S 2 in Chapter 7, p.~321 in [S]. Thus
$\Lambda_{c}^{\gamma}$ must only have one element. But this clearly proves the theorem since $c\in \Bbb
C^b$ is arbitrary. QED\vs We note in passing that we recover the theorem of A.A. Tarasev in [T] to
the effect that 
${\cal H}_y $ (see \S 3.3) is maximally Poisson commutative in $S(\g)$. Tarasev's result was in
response to the question of maximality posed by E. Vinberg.

\vs {\bf Theorem 4.12.} {\it The subalgebra
${\cal H}_y$ of \S 3.3 is maximally Poisson commutative in $S(\g)$.} 

\vs
\indent {\bf Proof.} Let $x\in \hbox{\rm Hess}$.
Then by the local Frobenius theorem (see Theorem 1.3.3, \S 1.3, Chapter 1, p.28 in [V] and the
statement of its applicabilty in the complex analytic case at the end of \S 1.3 on p.31) there exists
an complex open neighborhood $U'$ of $x$ in $\g^{\hbox{\sevenrm{sreg}}}$ which admits a foliation $$U' =
\sqcup_{\delta\in \Delta'} E_{\delta}\eqno (4.33)$$ where, for each $\delta$ in the parameter set
$\Delta'$, $E_{\delta}$ is a connected integral manifold for the distribution ${\cal Z}$
(see above in \S 4.2). Let $\Phi_U = \Phi|U'$ (see Theorem 4.1) and for any $z\in U'$ let
$\delta(z)\in \Delta'$ be such that $z\in E_{\delta(z)}$. Clearly $$E_{\delta(z)}\s F_c,\eqno (4.34)$$
where $c= \Phi_U(z)$. Now let $$U = \{z\in U'\mid \hbox{\rm Hess}\cap E_{\delta(z)} \neq \emptyset\}.$$ Note
that $U$ is not empty since $x\in U$. We assert that $U$ is complex open in $U'$. Indeed let $z\in U$
and let $c = \Phi_U(z)$. Then, recalling (4.9), one must have $v\in \hbox{\rm Hess}\cap E_{\delta(z)}$ where $$v =
(\Phi|\hbox{\rm Hess})^{-1}(c).\eqno (4.35)$$ But now $U'\cap \hbox{\rm Hess}$ is a complex open neighborhood of $v$ in 
$\hbox{\rm Hess}$.
By continuity there exists a complex open neighborhood $D_c$ of $c$ in $\Bbb C^b$ such that
$$(\Phi|\hbox{\rm Hess})^{-1}(D_c)\s U'\cap \hbox{\rm Hess}.$$ But by the continuity of $\Phi_U$ there exists a complex open
neighborhood $W$ of $z$ in $U'$ such that $\Phi_U(W)\s D_c$. But then it is immediate that $W\s U$.
Hence $U$ is open and clearly there exists a subset $\Delta\s \Delta'$ such that $$U =
\sqcup_{\delta\in \Delta} E_{\delta}.\eqno (4.36)$$

\vskip .1pc Now assume that $f\in S(\g)$ Poisson commutes with all $q\in {\cal H}_y$. But then by
Theorem 3.6 there exists $p\in {\cal H}_y$ such that $$f|\hbox{\rm Hess} = p| \hbox{\rm Hess}.$$ But both $p$ and $f$ are
constant on any connected integral manifold of ${\cal Z}$. But then $f-p$ vanishes on $U$ by (4.36).
Since $u$ is complex open in $\g$ this implies $f = p$. QED

\vs {\bf 4.3.} To state our final results it
will be convenient to replace $\Bbb C^b$ as parameters for the fibers $F_c$ of $\Phi|\g^{\hbox{\sevenrm{sreg}}}$ by
$\hbox{\rm Hess}$ (see (4.8) and (4.21)). For any $x\in \hbox{\rm Hess}$ put $$F_{[x]} = F_{\Phi(x)}.$$ We recall that
$$T_x(F_{[x]}) = Z_x;\eqno (4.37) $$ (see  (4.28) and Theorems 4.2 and 4.9).\vskip .5pc We also recall
that
$\hbox{\rm Hess}(O)$ is a Lagrangian submanifold of $O$ for any $O\in R$ (see Theorems 4.3 and 4.4). On the other
other hand one has 

\vs {\bf Theorem 4.13.} {\it Let $O\in R$. and let $x\in \hbox{\rm Hess}(O)$. Then
$F_{[x]}$ is a Lagrangian submanifold (not nessarily connected) of $O^{\hbox{\sevenrm{sreg}}}$. Furthermore $$T_x(O) = T
= T_x(F_{[x]}) \oplus T_x(\hbox{\rm Hess}(O))\eqno (4.38)$$ so that the two Lagrangian subspaces $T_x(F_{[x]})$
and 
$T_x(\hbox{\rm Hess}(O))$ of $T_x(O)$ are nonsingularly paired by $\omega_x$ (see (2.1))}.

 \vs {\bf Proof.} The
first conclusion of Theorem 4.13 is immediate from Theorems 4.2 and 4.9. To prove the final statement
of Theorem 4.13 it is, by dimension, enough to prove that $$T_x(F_{[x]}) \cap T_x(\hbox{\rm Hess}(O)) = 0.\eqno
(4.39) $$ But if $\Phi_*$ is the differential of $\Phi$ (operating on $T(\g)$), one has $T_x(F_{[x]})
\s \hbox{\rm Ker} \Phi_*$ by definition of the fiber $F_{[x]}$. On the other hand $\Phi_*|T_x(\hbox{\rm Hess}(O))$ is
injective by  Theorem 4.1. This proves (4.39).  QED

\vs For convenience, before our final statement we
will recount some of the definitions and previous results. We have chosen and fixed a regular
semisimple element $y\in
\hh$. Using $y$ and the translation of invariant procedure of Fomenko--Mischenko one constructs a maximal
Poisson commuatative subalgbra ${\cal H}_y$ of $S(\g)$. We use this Poisson commutative
subalgebra to introduce the definition of strong regularity in
$\g$ and the corresponding open Zariski dense variety $\g^{\hbox{\sevenrm{sreg}}}$ in $\g$. Intersecting with an adjoint
orbit $O$ of regular elements one has an open Zariski dense  subvariety $O^{\hbox{\sevenrm{sreg}}}$ of $O$
which then is also a symplectic submanifold of $O$. One has that
$\hbox{\rm Hess} \s \g^{\hbox{\sevenrm{sreg}}}$ where $\hbox{\rm Hess} = e_1 + \b_-$. Intersecting
  $\hbox{\rm Hess}$ with $O$ defines a
Lagrangian submanifold $\hbox{\rm Hess}(O)$ of $O^{\hbox{\sevenrm{sreg}}}$. The Hamiltonian vector fields which arise from
${\cal H}_y$, restricted to $\g^{\hbox{\sevenrm{sreg}}}$, define an involutive distribution ${\cal Z}$ on
$\g^{\hbox{\sevenrm{sreg}}}$. The leaves of ${\cal Z}$ define a foliation of $\g^{\hbox{\sevenrm{sreg}}}$.
  A choice of generators of
${\cal H}_y$ defines a surjective morphism $\Phi^{\hbox{\sevenrm{sreg}}}:\g^{\hbox{\sevenrm{sreg}}}\to \Bbb C^b$
  whose restriction to
$\hbox{\rm Hess}$ is an algbraic isomorphism $\hbox{\rm Hess} \to \Bbb C^b$. Here $b$ is the dimension of a Borel
  subalgebra
of $\g$. The irreducible components of any fiber $F_{[x]},\,x\in \hbox{\rm Hess}$,  of $\Phi^{\hbox{\sevenrm{sreg}}}$
  are maximal
complex connected integral submanifolds of ${\cal Z}$, establishing therefore that there are only a
finite number of maximal  complex connected integral submanifolds of ${\cal Z}$ in $F_{[x]}$ for any
$x\in \hbox{\rm Hess}$.

           The following result says in effect that then ${\cal H}_y$ simultaneously polarizes
$O^{\hbox{\sevenrm{sreg}}}$ for all regular orbits of $G$. Even more $O^{\hbox{\sevenrm{sreg}}}$ simulates a cotangent bundle structure
over a base manifold $\hbox{rm Hess}(O)$. \vs 

{\bf Theorem 4.14.} {\it Let $O\in R$
so that
$O$ is an arbitrary regular $G$-orbit in $\g$. Then $$O^{\hbox{\sevenrm{sreg}}} = \sqcup_{x\in \hbox{\rm Hess}(O)}
F_{[x]}\eqno (4.40)$$ defines a polarization of the symplectic manifold $O^{\hbox{\sevenrm{sreg}}}$, noting that
$\hbox{\rm Hess}(O)$ is a Lagrangian submanifold of $O^{\hbox{\sevenrm{sreg}}}$ and that $\hbox{\rm Hess}(O)$ is transversal to all
the Lagrangian fibers
$F_{[x]},\,x\in \hbox{\rm Hess}(O)$ of the polarization.}

\vs {\bf Proof.} Recalling the notation of (4.25) note
that $\Phi$ induces an isomorphism $\hbox{\rm Hess}(O) \to \Bbb C^b(O)$ by (a) of Theorem 4.3. But then (4.40)
follows from Theorem 4.8. The remaining statements follow from Theorem 4.4 and Theorem 4.13. QED

\vskip 8pt 
\centerline{\bf References}\vskip 4pt
\parindent=42pt
\item {[B]} A. Borel, {\it Linear Algebraic Groups}, W. A. Benjamin,  1969.
\item {[F-M]} A. Fomenko and A. Mishchenko, Euler equations on Lie
groups, {\it Math.~USSR-Izv.} {\bf 12} (1978), 371--389.
\item {[K1]} B. Kostant,  The Three-Dimensional Subgroup and the Betti Numbers of a Complex
 Simple Lie Group,  {\it Amer. J. Math.} {\bf 81}(1959), 973--1032.
\item {[K2]} B. Kostant, Lie Group Representations on Polynomial Rings, {\it Amer. J. Math.}, {\bf 85}(1963), 327--404.
\item {[K3]} B. Kostant, Quantization and Unitary
representations, {\it Lecture Notes in Math.}, Vol.~170, Springer-Verlag, 1970, pp.~87--207. 
\item{[K4]} Kostant, On Whittaker Vectors and Representation
Theory, Inventiones Math.,{\bf  48}(1978), 101--184.
\item {[M]} D. Mumford, {\it The Red Book of Varieties and Schemes},
Lecture Notes in Math., Vol.~1358, Springer, 1999.
\item {[S]} I. R. Shafarevitch, {\it Basic Algebraic Geometry},
 Springer, 1977.
\item {[Sp]} T. Springer, {\it Linear Algebraic Groups}, Second Edition.
 Progress in Math., Vol.~9, Birkh\"auser, 1998.
\item {[T]} A.A. Tarasev, {\it Russian Math. Surveys},{\bf 57} (2002),
  1013--1014
\item {[V]} V. S. Varadarajan, {\it Lie Groups, Lie Algebras and Their
Representations}, Prentice-Hall, 1974.

\smallskip
\parindent=30pt
\baselineskip=14pt
\vskip 1.9pc
\vbox to 60pt{\hbox{Bertram Kostant}
      \hbox{Dept. of Math.}
      \hbox{MIT}
      \hbox{Cambridge, MA 02139}
     \noindent E-mail:  kostant@math.mit.edu}

\rm

\end